\documentclass[11pt]{article}
\usepackage{natbib,enumerate}
 \bibpunct[, ]{(}{)}{,}{a}{,}{,}%

\usepackage{url}
\usepackage{amsmath,amsthm}    

\usepackage{amssymb}
\usepackage{fullpage}
\usepackage{fancyhdr}
\usepackage{setspace}
\usepackage[compact]{titlesec}
\newtheorem{theorem}{Theorem}

\newtheorem{lemma}{Lemma}
\newtheorem{proposition}{Proposition}
\newtheorem{remark}{Remark}

\newcommand{\ser}{s}
\newcommand{\sen}{n}
\newcommand{\del}{\delta_\lambda}
\newcommand{\delo}{\delta_\lambda}
\newcommand{\elo}{\eta_\lambda}
 
\newcommand{\rat}{\lambda/\mu} 

\begin{document}

\title{ Optimality Gap of Asymptotically-derived Prescriptions in Queueing Systems}
\author{Ramandeep S. Randhawa\\ \small{Marshall School of Business}\\ \small{University of Southern California}}
\date{ October 2012, Revised \today}
\maketitle

\begin{abstract}
In complex systems, it is quite common to resort to approximations when optimizing system performance. These approximations typically involve selecting a particular system parameter and then studying the performance of the system as this parameter grows without bound. In such an asymptotic regime, we prove that if the approximation to the objective function is accurate up to $\mathcal{O}(1)$, then under some regularity conditions,  the prescriptions that are derived from this approximation are $o(1)$-optimal, i.e., their optimality gap is asymptotically zero. A consequence of this result is that the well-known square-root staffing rules for capacity sizing in $M/M/s$  and $M/M/s+M$ queues to minimize the sum of linear expected steady-state customer waiting costs and linear capacity costs  are $o(1)$-optimal. We also discuss  extensions of this result for the case of non-linear customer waiting costs in these systems.
\end{abstract}

\section{Introduction}
\paragraph{Background.}
One typically encounters problems of optimizing  system capacity in queueing systems. With the exception of single server queueing systems such as $M/M/1$ or $M/G/1$, solving such problems in general is quite difficult. This difficulty arises due to the lack of a simple characterization of the performance measures of interest that is amenable to optimization. For instance, consider the problem of finding the optimal number of servers in an $M/M/s$ queueing system in steady-state when there is a customer-based holding cost of $h$ per customer per unit time and a capacity cost of $c$ per server per unit time; the goal being to minimize the sum of expected steady-state holding  and capacity costs.  That is, the objective is to solve
\begin{equation}
\label{eq:main-opt}
\min_{s\in \mathbb{Z}_+} \Pi(s):= h \mathbb{E}Q(s)+c s,
\end{equation}
where $\mathbb{E}Q(s)$ denotes the expected steady-state queue-length when there are $s$ servers, and $\mathbb{Z}_+$ denotes the set of non-negative integers. Denoting the customer arrival rate by $\lambda$ and the mean service time by $1/\mu$, it follows that the expected steady-state queue-length equals
\begin{equation}
\mathbb{E}Q(s)=\frac{\lambda}{s \mu-\lambda} B\left(\rho,s\right),
\end{equation}
where $\rho=\lambda/(s \mu)$ is the system utilization and $$B(\rho,s)=\frac{\frac{(s \rho)^{s}}{s!(1-\rho)}}{\frac{(s     \rho)^{s}}{s!(1-\rho)}+\sum_{k=0}^{s-1}\frac{(s \rho)^k}{k!}},~\text{for }s\in \mathbb{Z}_+$$ is the probability that all servers are busy. 

In solving \eqref{eq:main-opt},  we observe that the delay probability $B$ complicates the otherwise simple structure of the optimization problem and renders a straightforward solution impossible. To overcome this issue, the literature resorts to asymptotic analysis with the hope of obtaining an amicable approximation for $B$. (The details of this approach are in \cite{BMR:04}, but we reproduce some aspects for completeness.)  The asymptotic regime considered is that of large arrival rates, and  to make the dependence on the arrival rate $\lambda$ explicit, we denote the expected steady-state cost and queue-length by $\Pi_\lambda$ and $\mathbb{E}Q_\lambda$, respectively, and the  minimizer and optimal value function of \eqref{eq:main-opt} as $s_\lambda^\star$ and $\Pi_\lambda^\star$.

For our analysis, we use the following generalization of this delay probability formula that holds for all real $s>\lambda/\mu$ (see, for instance, \cite{JVD:86}):
\begin{equation}
  \label{eq:delay}
B(\rho,s)=\left[ s\rho\int_0^\infty  t e^{-s\rho t}(1+t)^{s-1}dt\right]^{-1}.
\end{equation}
We also note that for stability, we need the system capacity to exceed the arrival rate, i.e., $\mu s>\lambda$. It follows that we can limit attention to staffing levels of the form $s=\frac{\lambda}{\mu}+\hat{s}$, where $\hat{s}>0$ is typically referred to as the safety capacity, optimize \eqref{eq:main-opt} over the safety capacity $\hat{s}$.

As a first step in deriving the approximation-based  solution, one may wonder how many servers are needed as a proportion of the arrival rate. It is straightforward to establish that 
$$\lim_{\lambda\rightarrow \infty} \frac{s_\lambda^\star}{\lambda/\mu}=1,$$ that is, the optimal safety capacity as a fraction of the demand shrinks to zero. The linearity of the objective function in the capacity cost and expected steady-state queue-length further implies that the optimal number of servers must satisfy
$$
s_\lambda^\star=\frac{\lambda}{\mu}+\mathcal{O}(\sqrt{\lambda}).\footnote{For  functions   $f:\mathbb{R}\rightarrow \mathbb{R}$ and $g:\mathbb{R}\rightarrow \mathbb{R}_+$, we say that $f(\lambda)=\mathcal{O}(g(\lambda))$, if   $\limsup_{\lambda\rightarrow \infty} \frac{|f(\lambda)|}{g(\lambda)}<\infty$. Further, we say that   $f(\lambda)=o(g(\lambda))$ if $\lim_{\lambda\rightarrow \infty} \frac{f(\lambda)}{g(\lambda)}=0$.}
$$
In order to characterize the $\mathcal{O}(\sqrt{\lambda})$ term, a refined approximation is used. In particular, one focuses attention on number of servers of the form $s_\lambda=\frac{\lambda}{\mu}+ \sigma \sqrt{\frac{\lambda}{\mu} }$ and computes an approximation for the expected steady-state queue-length. For large arrival rates, the functional form of $s_\lambda$ chosen implies that the queue-length process behaves like a diffusion process, and the delay probability $B$ can be analyzed in more detail to yield
\begin{equation}
  \label{eq:Q-limit}
  \lim_{\lambda\rightarrow \infty} \frac{\mathbb{E}Q_\lambda\left(\frac{\lambda}{\mu}+ \sigma       \sqrt{\frac{\lambda}{\mu} }\right)}{\sqrt{\lambda/\mu}}=q(\sigma):=
\frac{1}{\sigma(1+ \sigma\Phi(\sigma)/\phi(\sigma))}  
\end{equation}
where $\phi$ and $\Phi$ denote the density and cumulative distribution function of the standard normal distribution (cf. \cite{HaW:81}). It follows that the ``best'' $\sigma$ can be found by solving
\begin{equation}
  \label{eq:approx-opt}
  \min_{\sigma>0} hq(\sigma)+ c \sigma.
\end{equation}
Using results in \cite{BMR:04}, it is easy to verify that the objective function above is strictly convex and hence noting that the objective function grows without bound as $\sigma$ approaches zero or infinity, we obtain that  \eqref{eq:approx-opt} has a unique solution.  Denoting this solution by $\sigma^\star$, we obtain a square-root staffing capacity prescription for the original problem
\begin{equation}
  \label{eq:ss}
  \bar{s}_\lambda=\frac{\lambda}{\mu}+\sigma^\star\sqrt{\frac{\lambda}{\mu}}.
\end{equation}
It follows that this diffusion-based approximation methodology leads us to within $o(\sqrt{\lambda})$ of the true optimal capacity (cf. \cite{JVZ:11}), i.e., we have
\begin{equation}
  \label{eq:ss-error}
  s_\lambda^\star=\bar{s}_\lambda+o(\sqrt{\lambda}).
\end{equation}
Further, given the asymptotic characterization in \eqref{eq:Q-limit}, we obtain that the optimality gap of the prescription $\bar{s}_\lambda$, defined as the difference between the costs incurred when using the prescription and when using the optimal value, is also $o(\sqrt{\lambda})$, i.e., we have
$$
\Pi_\lambda(\bar{s}_\lambda)=\Pi_\lambda^\star+o(\sqrt{\lambda}).
$$
The recent work \cite{JVZ:11} provides a further refined view of \eqref{eq:Q-limit} to prove that we actually have
\begin{equation}
  \label{eq:Q-limit-O(1)}
  \mathbb{E}Q_\lambda\left(\frac{\lambda}{\mu}+ \sigma \sqrt{\frac{\lambda}{\mu}     }\right)=q(\sigma)\sqrt{\frac{\lambda}{\mu}}+\mathcal{O}(1),
\end{equation}
which immediately yields that the square-root staffing $\bar{s}_\lambda$ is $\mathcal{O}(1)$-optimal, i.e.,
\begin{equation}
  \label{eq:O(1)-opt}
  \Pi_\lambda(\bar{s}_\lambda)=\Pi_\lambda^\star+\mathcal{O}(1).
\end{equation}
Note that $\bar{s}_\lambda$ may not be an integer, in which case the prescription could be chosen as either the closest integer smaller or larger than $\bar{s}_\lambda$ (depending on which performs better) and one expects \eqref{eq:O(1)-opt} to hold for the corresponding capacity level as well.

\paragraph{Main result of the paper.} In this paper, we prove that square-root staffing to minimize the sum of expected steady-state linear holding and capacity costs is in fact asymptotically exact or $o(1)$-optimal. That is, we have
\begin{equation}
  \label{eq:o(1)-opt}
  \lim_{\lambda\rightarrow \infty} \Big(\Pi_\lambda(\bar{s}_\lambda)-\Pi_\lambda^\star\Big)=0.
\end{equation}
We in fact prove the following stronger result that allows small deviations from the square-root staffing prescription:
\begin{equation}
  \label{eq:o(1)-opt-general}
  \lim_{\lambda\rightarrow \infty} \Big(\Pi_\lambda(\bar{s}_\lambda+\eta_\lambda)-\Pi_\lambda^\star\Big)=0,~\text{ for }\eta_\lambda=o(\lambda^{1/4}).
\end{equation}
Thus, as the arrival rate $\lambda$ grows, the optimal cost grows along with it, however, surprisingly, the optimality gap of the square-root staffing prescription decreases and is asymptotically zero. To the best of my knowledge, this paper is the first to establish that the unscaled optimality gap of a capacity prescription derived using asymptotic analysis in a queueing system diminishes to zero.

The intuition behind this extremely accurate performance stems from the fact that objective functions tend to be flat in the vicinity of their optimizers. However, this intuition needs to be refined in the face of the $\mathcal{O}(1)$ error term in \eqref{eq:Q-limit-O(1)} because if this term is not well behaved it can lead to an overall $\mathcal{O}(1)$ optimality gap. One way to establish that this error term is well behaved is to explicitly compute it and then analyze it. For instance, \cite{JVZ:11} does compute this term, and this can be used to establish the result. However, this computation is indeed quite intricate and it is unclear if it is doable for general systems. In this paper, we utilize properties of the underlying queueing system to establish that if the error term is $\mathcal{O}(1)$-accurate, then the prescription is $o(1)$-optimal.

\paragraph{Other results}
Though our focus is on square-root staffing to minimize linear costs in $M/M/s$ systems in steady-state, our methods easily extend to other asymptotic optimization problems as long as the performance measure of interest has some regularity properties and has an $\mathcal{O}(1)$-accurate approximation (the recent work \cite{Gurvich:14}  proves that such an approximation result holds for fairly general Markovian systems). In particular, we prove that square-root staffing to minimize linear costs is also $o(1)$-optimal for $M/M/s+M$ systems in steady-state and that the fluid-approximation based prescriptions derived in \cite{BaR:10}  are $o(1)$-optimal for overloaded $M/M/s+GI$ systems in steady-state.

We also consider the case of non-linear delay costs of the power form, i.e., a customer waiting for $W$ time units leads to a system cost of the form $h W^r$ for $r>0$. For $M/M/s+M$ systems with convex delay costs, we obtain somewhat surprisingly that fluid approximations are $o(1)$-optimal. For $M/M/s$ systems, the result is more nuanced because the safety capacity needed in excess of the demand, depends on $r$ and further, so does the approximation error. We characterize this approximation error as $\mathcal{O}(n_\lambda)$, where $n_\lambda=\lambda^{\frac{3-r}{2(1+r)}}$, and our results imply that the corresponding prescription has an optimality gap that is an order smaller;  in particular the prescription is $o(n_\lambda)$-optimal.

\section{Square-root staffing is asymptotically exact}
\subsection{The result}
In this section, we focus on the cost minimization problem \eqref{eq:main-opt} in the $M/M/s$ system and prove that the corresponding square-root staffing rule has asymptotically exact performance. We will use the following asymptotic property for $M/M/s$ queueing systems that follows from Theorem~2 in \cite{JVZ:11}:
\begin{lemma}
  \label{lem:mms}
  For any $\sigma>0$, we have
  \begin{equation}
    \label{eq:ss-ass}
    \mathbb{E}Q_\lambda\left(\frac{\lambda}{\mu}+\sigma       \sqrt{\frac{\lambda}{\mu}}\right)=q(\sigma)\sqrt{\frac{\lambda}{\mu}}+\epsilon_\lambda(\sigma),
  \end{equation}
  where $q:\mathbb{R}_+\rightarrow \mathbb{R}_+$ (as defined in \eqref{eq:Q-limit}), $\epsilon_\lambda:\mathbb{R}_+\rightarrow   \mathbb{R}$ and:
\begin{enumerate}[(a)]
\item $q$ is twice continuously differentiable.
\item $\epsilon_\lambda(\sigma)=o(\sqrt{\lambda})$, i.e., the sequence $\{\epsilon_\lambda(\sigma)\}$ satisfies  $ \lim_{\lambda\rightarrow \infty} \frac{\epsilon_\lambda(\sigma)}{\sqrt{\lambda}}=0$.
\item \label{item:imp} (Continuity of error term) For any real-valued sequence $\{\kappa_\lambda\}$ such that $\lim_{\lambda \rightarrow \infty} \kappa_\lambda=0$, we have  $\lim_{\lambda \rightarrow \infty} \big(\epsilon_\lambda(\sigma^\star+\kappa_\lambda)-\epsilon_\lambda(\sigma^\star)\big)=0$.
\end{enumerate}
\end{lemma}
Note that establishing property~\eqref{item:imp}  is the most difficult part of the analysis because it requires analyzing the error term $\epsilon_\lambda$ and proving that it does not change very rapidly. We next use Lemma~\ref{lem:mms} to establish the asymptotic optimality property of square-root staffing and then we will return to analyzing property~\eqref{item:imp}  in more detail.
\begin{theorem}
\label{thm:main}
Any staffing level  $\bar{s}_\lambda+\eta_\lambda$, where $\bar{s}_\lambda$ is the square-root staffing prescription defined in \eqref{eq:ss} and $\eta_\lambda=o(\lambda^{1/4})$, is $o(1)$-optimal for the optimization problem ~\eqref{eq:main-opt}.
\end{theorem}
\begin{proof}
We first establish the following result that square-root staffing is $o(1)$-optimal:
\begin{equation}
\label{eq:liminf-0}
\lim_{\lambda \rightarrow \infty} \Big(\Pi_\lambda(\bar{s}_\lambda)-\Pi_\lambda^\star\Big)=0.
\end{equation}
Using \eqref{eq:ss-error}, we can write $s_\lambda^\star=\bar{s}_\lambda+\delo$ for some $\delo=o(\sqrt{\lambda})$.
Defining $\hat{\pi}(\sigma)=h q(\sigma)+c \sigma$, we have
\begin{align}
\nonumber
\Pi_\lambda^\star&=\Pi_\lambda(\bar{s}_\lambda+\delo)\\
&=h \mathbb{E}Q_\lambda(\bar{s}_\lambda+\delo)+c (\bar{s}_\lambda+\delo)\\
\nonumber
&=c \frac{\lambda}{\mu}+h \mathbb{E}Q_\lambda\left(\frac{\lambda}{\mu}+ (\sigma^\star+\frac{\delo}{\sqrt{\lambda/\mu}})\sqrt{\frac{\lambda}{\mu} }\right)+ c \sqrt{\frac{\lambda}{\mu}}\left(\sigma^\star+\frac{\delo}{\sqrt{\lambda/\mu}}\right)\\
\nonumber
&=c \frac{\lambda}{\mu}+ \sqrt{\frac{\lambda}{\mu} }\left[h q\left(\sigma^\star+\frac{\delo}{\sqrt{\lambda/\mu}}\right)+c\left(\sigma^\star+\frac{\delo}{\sqrt{\lambda/\mu}}\right)\right]+h \epsilon_\lambda\left(\sigma^\star+\frac{\delo}{\sqrt{\lambda/\mu}}\right)\\
\label{eq:optimal}
&=c \frac{\lambda}{\mu} +\sqrt{\frac{\lambda}{\mu} }\hat{\pi}\left(\sigma^\star+\frac{\delo}{\sqrt{\lambda/\mu}}\right)+h \epsilon_\lambda\left(\sigma^\star+\frac{\delo}{\sqrt{\lambda/\mu}}\right).
\end{align}
We also have
\begin{equation}
\label{eq:ss-cost}
\Pi_\lambda(\bar{s}_\lambda)=c\frac{\lambda}{\mu}+\sqrt{\frac{\lambda}{\mu} }\hat{\pi}(\sigma^\star)+h \epsilon_\lambda(\sigma^\star).
\end{equation}
Thus, comparing \eqref{eq:optimal} and \eqref{eq:ss-cost}, we have
\begin{align}
\nonumber
\Pi_\lambda(\bar{s}_\lambda)-\Pi_\lambda^\star=&\sqrt{\frac{\lambda}{\mu} }\left[\hat{\pi}(\sigma^\star)-\hat{\pi}\left(\sigma^\star+\frac{\delo}{\sqrt{\lambda/\mu}}\right)\right]\\
\label{eq:subs}
&+h \left[\epsilon_\lambda(\sigma^\star)-\epsilon_\lambda\left(\sigma^\star+\frac{\delo}{\sqrt{\lambda/\mu}}\right)\right].
\end{align}

Consider the first term in \eqref{eq:subs}. Using the fact that $\sigma^\star$ minimizes $\hat{\pi}$, we have
\begin{equation}
\label{eq:term1}
\sqrt{\frac{\lambda}{\mu} }\left[\hat{\pi}(\sigma^\star)-\hat{\pi}\left(\sigma^\star+\frac{\del}{\sqrt{\lambda/\mu}}\right)\right]\le 0.
\end{equation}
Turning to the second term in \eqref{eq:subs} and applying Lemma~\ref{lem:mms}\eqref{item:imp}, we obtain that
\begin{equation}
\label{eq:term2}
\lim_{\lambda\rightarrow \infty} \left[\epsilon_\lambda(\sigma^\star)-\epsilon_\lambda\left(\sigma^\star+\frac{\del}{\sqrt{\lambda/\mu}}\right)\right]=0.
\end{equation}
Thus, combining \eqref{eq:term1} and \eqref{eq:term2} in \eqref{eq:subs}, we obtain
\begin{equation}
\label{eq:liminf}
\limsup_{\lambda \rightarrow \infty} \Big(\Pi_\lambda(\bar{s}_\lambda)-\Pi_\lambda^\star\Big)\le 0.
\end{equation}
Because $\Pi_\lambda^\star$ is the optimal cost, we have  $\Pi_\lambda(\bar{s}_\lambda)-\Pi_\lambda^\star\ge 0$, and thus, \eqref{eq:liminf-0} holds, and we obtain that square-root staffing is $o(1)$-optimal.

We next establish the $o(1)$-optimality for the staffing level $\bar{s}_\lambda+\eta_\lambda$. To do so, we will prove that 
$$
\lim_{\lambda \rightarrow \infty}\left(  \Pi_\lambda(\bar{s}_\lambda)-\Pi_\lambda(\bar{s}_\lambda+\eta_\lambda)\right)=0.$$
We will use the arguments used previously in this proof for the staffing level $\bar{s}_\lambda$ until \eqref{eq:subs} because these only use  $\Pi_\lambda^\star=\Pi_\lambda(\bar{s}_\lambda+\delta_\lambda)$ for some $\delta_\lambda=o(\sqrt{\lambda})$ and do not use the optimality of $\Pi_\lambda^\star$. This gives us the following analog of \eqref{eq:subs}:
\begin{align}
\nonumber
\Pi_\lambda(\bar{s}_\lambda)-\Pi_\lambda(\bar{s}_\lambda+\eta_\lambda)=&\sqrt{\frac{\lambda}{\mu} }\left[\hat{\pi}(\sigma^\star)-\hat{\pi}\left(\sigma^\star+\frac{\eta_\lambda}{\sqrt{\lambda/\mu}}\right)\right]\\
\label{eq:subs2}
&+h \left[\epsilon_\lambda(\sigma^\star)-\epsilon_\lambda\left(\sigma^\star+\frac{\elo}{\sqrt{\lambda/\mu}}\right)\right].
\end{align}
As in \eqref{eq:term2}, the second term above converges to zero as $\lambda$ grows without bound. We next apply the Taylor series expansion to the first term in \eqref{eq:subs2}. This yields
\begin{equation*}
\sqrt{\frac{\lambda}{\mu} }\left[\hat{\pi}(\sigma^\star)-\hat{\pi}\left(\sigma^\star+\frac{\elo}{\sqrt{\lambda/\mu}}\right)\right]
=\sqrt{\frac{\lambda}{\mu} }\left[-\hat{\pi}'(\sigma^\star) \frac{\elo}{\sqrt{\lambda/\mu}}-\hat{\pi}''(\xi_\lambda)\frac{1}{2}\left(\frac{\elo}{\sqrt{\lambda/\mu}}\right)^2\right],
\end{equation*}
where $\xi_\lambda \in (\sigma^\star,\sigma^\star+{\elo}/{\sqrt{\lambda/\mu}})$. Noting that $\sigma^\star$ minimizes $\hat{\pi}$ and that $\hat{\pi}(x)\rightarrow \infty$ as $x\rightarrow 0,\infty$, we must have $\hat{\pi}'(\sigma^\star)=0$. This gives us
\begin{equation}
\label{eq:term1-2}
\sqrt{\frac{\lambda}{\mu} }\left[\hat{\pi}(\sigma^\star)-\hat{\pi}\left(\sigma^\star+\frac{\elo}{\sqrt{\lambda/\mu}}\right)\right]
= -\hat{\pi}''(\xi_\lambda) \frac{\eta_\lambda^2}{2\sqrt{\lambda/\mu}}.
\end{equation}
Thus taking limits as $\lambda\rightarrow \infty$ in \eqref{eq:subs2} and using the fact that $\hat{\pi}(\sigma)=h q(\sigma)+c \sigma$ is twice continuously differentiable at $\sigma^\star$ and that $\hat{\pi}''(\sigma^\star)\ge 0$, we obtain
$$\liminf_{\lambda \rightarrow \infty}\left(  \Pi_\lambda(\bar{s}_\lambda)-\Pi_\lambda(\bar{s}_\lambda+\eta_\lambda)\right)=-\hat{\pi}''(\sigma^\star)\limsup_{\lambda\rightarrow \infty} \frac{\eta_\lambda^2}{2 \sqrt{\lambda/\mu}}.
$$
Because $\eta_\lambda=o(\lambda^{1/4})$, we obtain that 
$$\liminf_{\lambda \rightarrow \infty}\left(  \Pi_\lambda(\bar{s}_\lambda)-\Pi_\lambda(\bar{s}_\lambda+\eta_\lambda)\right)=0.
$$
A similar argument holds using ``lim sup'' instead of ``lim inf'' so that we obtain 
$$\limsup_{\lambda \rightarrow \infty}\left(  \Pi_\lambda(\bar{s}_\lambda)-\Pi_\lambda(\bar{s}_\lambda+\eta_\lambda)\right)=0.
$$
Thus we have $$
\lim_{\lambda \rightarrow \infty}\left(  \Pi_\lambda(\bar{s}_\lambda)-\Pi_\lambda(\bar{s}_\lambda+\eta_\lambda)\right)=\lim_{\lambda \rightarrow \infty}\left(  \Pi_\lambda(\bar{s}_\lambda)-\Pi_\lambda^\star\right)=0.
$$
$\qed$\end{proof}

Theorem~\ref{thm:main} proves that square-root staffing is asymptotically $o(1)$-optimal. In fact, \cite{JVZ:11} proves that $\epsilon_\lambda(\sigma)=\epsilon(\sigma)+\mathcal{O}(1/\sqrt{\lambda})$, for some function $\epsilon$. Using this in the proof of Theorem~\ref{thm:main}, we can establish the stronger result that the optimality gap of square-root staffing is $\mathcal{O}(1/\sqrt{\lambda})$.

\subsection{Drivers of asymptotic performance}
Theorem~\ref{thm:main} proves that even though for any square-root staffing level, the diffusion approximation for the system cost is  within $\mathcal{O}(1)$ of the actual cost, the cost under the  square-root staffing prescription is within $o(1)$ of the optimal value. The goal in this section is to better understand the driver for the fact that a prescription obtained from asymptotic methods performs an `order' better than what one would naively expect.

The key property in Lemma~\ref{lem:mms}  used to establish Theorem~\ref{thm:main} is the  (limiting) continuity of the error term, in particular, that for $\del=\mathcal{O}({\lambda}^{1/4})$, we have
\begin{equation}
\label{eq:relation}
\lim_{\lambda\rightarrow \infty} \left[\epsilon_\lambda(\sigma^\star)-\epsilon_\lambda\left(\sigma^\star+\frac{\del}{\sqrt{\lambda/\mu}}\right)\right]=0.
\end{equation}
Notice that this continuity result is only needed for $ \del=\mathcal{O}(\lambda^{1/4})$ rather than for all $\del=o(\sqrt{\lambda})$. This is so because \eqref{eq:subs} and \eqref{eq:term1-2} imply that
\begin{align*}
0\le \Big( \Pi_\lambda(\bar{s}_\lambda)-\Pi_\lambda^\star\Big)=& -\hat{\pi}''(\xi_\lambda)\frac{\del^2}{2\sqrt{\lambda/\mu}}
+h\left[\epsilon_\lambda(\sigma^\star)-\epsilon_\lambda\left(\sigma^\star+\frac{\del}{\sqrt{\lambda/\mu}}\right)\right],
\end{align*}
and noting that $\pi''(\xi_\lambda)\rightarrow \pi''(\sigma^\star)>0$ and that $\epsilon_\lambda$ is bounded, the term $\del^2/\sqrt{\lambda}$ must also be bounded (otherwise the optimality of $\Pi_\lambda^\star$ would be contradicted). 

Lemma~\ref{lem:mms} was established using  results in \cite{JVZ:11}, which proved that $\epsilon_\lambda=\mathcal{O}(1)$ and also characterized this term explicitly. We will show how \eqref{eq:relation} can be directly established for queueing systems by  using the property $\epsilon_\lambda=\mathcal{O}(1)$, without any additional characterization of the error term, and using another basic property of the expected steady-state queue-length. In particular, we  will show in Proposition~\ref{prop:intuition}, that the curvature of the  expected steady-state queue length as a function of the number of servers can be used to establish this continuity property. Intuitively,  for small changes to the number of servers, the curvature of the expected steady-state queue-length is ``picked up'' solely by the error term $\epsilon_\lambda$, and hence the properties of the curvature can be used to prove that for these small changes $\epsilon_\lambda$ must be well-behaved. In fact, the well-known convexity of the expected steady-state queue-length can be used to establish \eqref{eq:relation} for $\del=o(\lambda^{1/4})$. To extend this argument and cover the case of all  $\del=\mathcal{O}(\lambda^{1/4})$, we need additional properties about the curvature of the queue-length. In particular, Proposition~\ref{prop:mmn-convex} proves that the convexity of the expected steady-state queue-length is increasing with utilization (for fixed number of servers)  and decreasing in the number of servers (for any fixed offered load, $\lambda/\mu$). This additional property, that we state below, allows us to establish \eqref{eq:relation}  for $\del=\mathcal{O}(\lambda^{1/4})$. 
\begin{proposition}
\label{prop:mmn-convex}
For an $M/M/s$ system, for fixed number of servers $s$, the convexity of the expected steady-state queue-length with respect to the system utilization $\rho=\frac{\lambda}{s \mu}$ is increasing in the system utilization. Further, for any fixed offered load $\lambda/\mu$, the convexity of the expected steady-state queue-length with respect to the number of servers $s$ is decreasing with the number of servers, for $s\ge 3$.
\end{proposition}
The proof of this result is straightforward, though lengthy, and is postponed to the Appendix. We next show how Proposition~\ref{prop:mmn-convex} allows us to establish \eqref{eq:relation}.
\begin{proposition}
\label{prop:intuition}
For an $M/M/s$ queueing system, if the approximation \eqref{eq:ss-ass} holds for $q$ that is thrice continuously differentiable and $\epsilon_\lambda(\sigma)=\mathcal{O}(1)$ for all $\sigma>0$, then for any $\sigma>0$ and $\del=\mathcal{O}(\lambda^{1/4})$ we have
\begin{equation}
\label{eq:continuity}
\lim_{\lambda\rightarrow \infty}\left|\epsilon_\lambda(\sigma)-\epsilon_\lambda\left(\sigma+\frac{\del}{\sqrt{\lambda/\mu}}\right)\right|=0.
\end{equation}

\end{proposition}
\begin{proof}
Fix any $\sigma>0$ and sequence $\del=
\mathcal{O}(\lambda^{1/4})$. We begin with some definitions. For any integer $i$ such that $\lambda/\mu+\sigma \sqrt{\lambda/\mu}+i\del\ge 0$, we define
$$\Delta_\lambda^\epsilon(i):=\epsilon_\lambda\left(\sigma+(i+1)\frac{\del}{\sqrt{\lambda/\mu}}\right)-\epsilon_\lambda\left(\sigma+i \frac{\del}{\sqrt{\lambda/\mu}}\right).
$$
Notice that because $\epsilon_\lambda$ is bounded,  the sequence $\Delta_\lambda^\epsilon(i)$ for any fixed $i$ is also bounded, and we can define
$$
\alpha:=\limsup_{\lambda\rightarrow \infty}\Delta_\lambda^\epsilon(0).
$$
For ease of notation, we can assume that  $\lim_{\lambda\rightarrow \infty} \Delta_\lambda^\epsilon(0)=\alpha$ (otherwise we can proceed with the corresponding convergent subsequence). To complete the proof, we need to prove that $\alpha=0$. 
We next define 
$$
K:=\sup_{x\in[\sigma/2,2\sigma]}\limsup_{\lambda\rightarrow \infty} |\epsilon_\lambda(x)|,
$$
where we take the supremum in a neighborhood of $\sigma$ to ensure that the sequence $\sigma +i \frac{\del}{\sqrt{\lambda/\mu}}$ is covered in the neighborhood asymptotically. The definition of $K$ ensures that we have $\limsup_{\lambda\rightarrow \infty} |\Delta_\lambda^\epsilon(i)|\le 2K$ for all $i$. 

We next compute the difference in expected steady-state queue-length obtained by a small change in the number of servers as follows:
\begin{align}
\Delta_\lambda^Q(i):=&\mathbb{E}Q_\lambda\left(\frac{\lambda}{\mu}+\sigma \sqrt{\frac{\lambda}{\mu}}+(i+1) \del\right)-
\mathbb{E}Q_\lambda\left(\frac{\lambda}{\mu}+\sigma \sqrt{\frac{\lambda}{\mu}}+i  \del\right)\nonumber\\
&=\sqrt{\frac{\lambda}{\mu}} q\left(\sigma+(i+1)\frac{\del}{\sqrt{\lambda/\mu}}\right)-\sqrt{\frac{\lambda}{\mu}} q\left(\sigma+i \frac{\del}{\sqrt{\lambda/\mu}}\right)+\Delta_\lambda^\epsilon(i)\nonumber\\
&\stackrel{(a)}{=}q'(\sigma) \del+(2 i+1)q''(\sigma) \frac{\del^2}{2\sqrt{\lambda/\mu}}
+\big((i+1)^3 q'''(\psi_\lambda)-i^3 q'''(\xi_\lambda)\big)\frac{\del^3}{6{\lambda/\mu}}\nonumber\\
&+\Delta_\lambda^\epsilon(i),\nonumber\\
&\stackrel{(b)}{=}q'(\sigma) \del+(2 i+1)q''(\sigma) \frac{\del^2}{2\sqrt{\lambda/\mu}}+\big((i+1)^3-i^3\big) q'''(\sigma)\frac{\del^3}{6{\lambda/\mu}}\nonumber\\
&+o(\lambda^{-1/4})+\Delta_\lambda^\epsilon(i),
\label{eq:DQ1}
\end{align}
where, in $(a)$,  $\psi_\lambda$ lies between $\sigma$ and $\sigma+(i+1){\del}/{\sqrt{\lambda/\mu}}$ and  $\xi_\lambda$ lies between $\sigma$ and $\sigma+i{\del}/{\sqrt{\lambda/\mu}}$, and in $(b)$, we use that $\del=\mathcal{O}(\lambda^{1/4})$ and that $q'''$ is continuous. Notice that $\Delta_\lambda^Q$ is related to the negative of the first derivative of $\mathbb{E}Q_\lambda(s)$. Using the expression for $\Delta_\lambda^Q(i)$ from \eqref{eq:DQ1}, we obtain the following relation, which relates to the second derivative of the expected queue-length:
\begin{align}
\Delta_\lambda^Q(i+1)-\Delta_\lambda^Q(i)=&q''(\sigma) \frac{\del^2}{\sqrt{\lambda/\mu}}+(i+1)q'''(\sigma)\frac{\del^3}{{\lambda/\mu}}\nonumber
\\
&+\Delta_\lambda^\epsilon(i+1)-\Delta_\lambda^\epsilon(i)+o(\lambda^{-1/4}).\label{eq:DQ}
\end{align}

We next utilize the result from Proposition~\ref{prop:mmn-convex} that the convexity (second-derivative) of the expected steady-state queue-length with respect to the number of servers is decreasing. This gives us the following relation:
\begin{align}
\label{eq:decreasing-DQ}
\Delta_\lambda^Q(i+2)-\Delta_\lambda^Q(i+1)&\le \Delta_\lambda^Q(i+1)-\Delta_\lambda^Q(i).
\end{align}
We next use \eqref{eq:DQ} for the right-hand-side of the above relation, and the following relation for the left-hand-side that is obtained analogous to \eqref{eq:DQ} by replacing $(i+1), i$ by $(i+2), (i+1)$ respectively:
\begin{align}
\Delta_\lambda^Q(i+2)-\Delta_\lambda^Q(i+1)=&q''(\sigma) \frac{\del^2}{\sqrt{\lambda/\mu}}+(i+2)q'''(\sigma)\frac{\del^3}{{\lambda/\mu}}\nonumber
\\
&+\Delta_\lambda^\epsilon(i+2)-\Delta_\lambda^\epsilon(i+1)+o(\lambda^{-1/4}).\label{eq:DQ+1}
\end{align}
This gives us the following relation:
\begin{align}
\Delta_\lambda^\epsilon(i+2)-\Delta_\lambda^\epsilon(i+1)&\le \Delta_\lambda^\epsilon(i+1)-\Delta_\lambda^\epsilon(i)+o(\lambda^{-1/4})- q'''(\sigma)\frac{\del^3}{{\lambda/\mu}}.
\end{align}

For $\del=\mathcal{O}(\lambda^{1/4})$, thus it follows that
\begin{equation}
\label{eq:third-D}
\limsup_{\lambda\rightarrow \infty} \Big(\big[\Delta_\lambda^\epsilon(i+2)-\Delta_\lambda^\epsilon(i+1)\big]-\big[ \Delta_\lambda^\epsilon(i+1)-\Delta_\lambda^\epsilon(i)\big]\Big)\le 0.
\end{equation}

We will next prove that
$$
\lim_{\lambda\rightarrow \infty} \left(\Delta_\lambda^\epsilon(i+1)-\Delta_\lambda^\epsilon(i)\right)=0 \text{ for all }i. 
$$

Toward a contradiction, suppose that $\liminf_{\lambda\rightarrow \infty} \left(\Delta_\lambda^\epsilon(i+1)-\Delta_\lambda^\epsilon(i)\right)=\beta<0$ for some integer $i$ (where $\beta$ may depend on $i$). Then, we use \eqref{eq:third-D} to obtain 
\begin{equation*}
\liminf_{\lambda\rightarrow \infty} \left( \Delta_\lambda^\epsilon(m+i)-\Delta_\lambda^\epsilon(i) \right)\le m \beta \text{ for } m=1,2,\dots.
\end{equation*}

Thus, for $m>|4 K/\beta|$, we obtain $\liminf_{\lambda\rightarrow \infty} |\Delta_\lambda^\epsilon(m+i)-\Delta_\lambda^\epsilon(i)|>4 K$, which contradicts the fact that $|\Delta_\lambda^\epsilon(m+i)-\Delta_\lambda^\epsilon(i)|\le 4K$, which follows because $|\Delta_\lambda^\epsilon(\ell)|\le  2 K$ for all $\ell$. Thus, we must have  
\begin{equation*}
\liminf_{\lambda\rightarrow \infty} \left(\Delta_\lambda^\epsilon(i+1)-\Delta_\lambda^\epsilon(i)\right)\ge 0 \text{  for all } i.
\end{equation*}

Next, suppose that we have $\limsup_{\lambda\rightarrow \infty} \left(\Delta_\lambda^\epsilon(i+1)-\Delta_\lambda^\epsilon(i)\right)=\beta>0$ for some integer $i$. Then, proceeding as in the previous argument, we obtain that $\limsup_{\lambda\rightarrow \infty} \left(\Delta_\lambda^\epsilon(i+1)-\Delta_\lambda^\epsilon(i+1-m)\right)\ge m \beta$ for $m=1,2,\dots$, and a corresponding contradiction is obtained. Thus, we can conclude that  
\begin{equation*}
\limsup_{\lambda\rightarrow \infty} \left(\Delta_\lambda^\epsilon(i+1)-\Delta_\lambda^\epsilon(i)\right)\le 0 \text{ for all }i.
\end{equation*}
This in fact proves that we must have \begin{equation*}
\lim_{\lambda\rightarrow \infty} \left(\Delta_\lambda^\epsilon(i+1)-\Delta_\lambda^\epsilon(i)\right)=0 \text{ for all }i. 
\end{equation*}
This also gives us $\lim_{\lambda \rightarrow \infty} \Delta_\lambda^\epsilon(i)=\lim_{\lambda \rightarrow \infty} \Delta_\lambda^\epsilon(0)=\alpha$ for all $i$. Thus, $$\lim_{\lambda \rightarrow \infty} \left[\epsilon_\lambda\left(\sigma+m\frac{\del}{\sqrt{\lambda/\mu}}\right)-\epsilon_\lambda\left(\sigma \right)\right]=\lim_{\lambda \rightarrow \infty}  \sum_{i=0}^{m-1} \Delta_\lambda^\epsilon(i)=m \alpha.$$ 
Using the definition of $K$, it follows that we must have $|m \alpha|\le 2 K$ for all $m$. This implies that we must have $\alpha=0$, which completes the proof.
$\qed$\end{proof}

\begin{remark}[The role of Proposition~\ref{prop:mmn-convex}.]
Proposition~\ref{prop:mmn-convex} establishes  that the convexity of the expected queue-length with respect to the number of servers is decreasing. We would  like to point out that although establishing this property required some work and is new to the literature, the property is not surprising and was expected. In this sense,  we expect the  order improvement in performance of asymptotically derived prescriptions that we observe here should apply to other systems as well.

Technically, we would also like to point out that the proof in Proposition~\ref{prop:intuition} works under weaker conditions on the convexity. In particular,  we only need one of the following conditions to hold for each $\sigma>0$: 
\begin{align}
\label{eq:third-weak}
\liminf_{\lambda\rightarrow \infty} D^{(3)}Q_\lambda\left(\frac{\lambda}{\mu}+\sigma \sqrt{\frac{\lambda}{\mu}}+\del \right)&\ge 0, \text{ for all } \delta_\lambda=\mathcal{O}(\lambda^{1/4}),\text{ or }\\
\limsup_{\lambda\rightarrow \infty} D^{(3)}Q_\lambda\left(\frac{\lambda}{\mu}+\sigma \sqrt{\frac{\lambda}{\mu}}+\del \right)&\le 0,  \text{ for all } \delta_\lambda=\mathcal{O}(\lambda^{1/4}),
\end{align}
where $D^{(j)} Q_\lambda(s)=D^{(j-1)} Q_\lambda(s+1)-D^{(j-1)} Q_\lambda(s),$ for $j=1,2,3$ 
with $D^{(0)} Q_\lambda(s)=\mathbb{E} Q_\lambda(s)$. That is,  proving Proposition~\ref{prop:intuition} only requires the third derivative to not change signs for $\mathcal{O}(\lambda^{1/4})$ changes in the number of servers.
\end{remark}

\section{Extensions}
In this section, we discuss how our results extend to non-linear delay costs in  $M/M/s$ queueing systems and to systems with customer abandonments for both linear and non-linear delay costs. For convenience, we will recycle some notation from the previous section. In particular, in each setting,  $\bar{s}_\lambda$ denotes the capacity prescription, $q(\sigma)$ denotes the approximation to the expected steady-state queue-length, and $\sigma^\star$ denotes the solution to the approximate optimization problem.

\subsection{Non-linear delay costs in $M/M/s$ systems}
\label{sec:non-linear}
We will focus on the following non-linear version of ~\eqref{eq:main-opt} that was also studied in \cite{KuR:10}:
\begin{equation}
\label{eq:main-opt-nl}
\min_{s\in \mathbb{Z}_+} \Pi_\lambda(s):= h \lambda \mathbb{E}\xi_\lambda(s)+c s,
\end{equation}
where $\mathbb{E}\xi_\lambda(s)=\mathbb{E}W_\lambda(s)^r$ denotes the $r^{th}$ moment of the steady-state waiting time or time-in-queue that the customers experience. Using the steady-state waiting time distribution for an $M/M/s$ queue (cf. \cite{CL:03}), we have
\begin{equation}
  \label{eq:1}
   \mathbb{E}W_\lambda(s)^r=\frac{\Gamma(r+1)}{(s\mu-\lambda)^r}B(\rho,s).
\end{equation}
As analyzed in depth in \cite{KuR:10}, for $r\ne 1$, the optimal staffing solution for \eqref{eq:main-opt-nl} is not square-root staffing. As we will soon see, the optimality gaps here are also different compared with the linear case. We analyze the case of convex delay costs ($r>1$) in detail in Section~\ref{sec:convex} and then discuss the case of concave delay costs ($r<1$) in Section~\ref{sec:concave}.

\subsubsection{Convex delay costs: $r>1$.}
\label{sec:convex} In this case, the optimal staffing level is smaller than the square-root level and is given by
\begin{equation}
  \label{eq:staff-nl}
  {s}^\star_\lambda=\frac{\lambda}{\mu}+{\sigma}^\star\left( \frac{\lambda}{\mu} \right)^{\frac{1}{r+1}}+\epsilon_\lambda,
\end{equation}
where
$$
{\sigma}^\star=\left( \frac{h}{c} r\Gamma(r+1) \mu^{1-r}\right)^{\frac{1}{r+1}},
$$
and $\epsilon_\lambda=o\left(\lambda^{\frac{1}{r+1}}\right)$ is the error term. 
Intuitively for $r>1$, if we use square-root staffing, i.e., $s_\lambda=\frac{\lambda}{\mu}+\mathcal{O}(\sqrt{\lambda})$, then the waiting costs equal $\mathcal{O}(\lambda^{1-r/2})=o(\sqrt{\lambda})$, and thus are an order smaller than the (safety) capacity costs. This implies that to optimize the overall cost, capacity needs to be reduced to a smaller order. Further, with a reduced size of the safety capacity (any order smaller than that in square-root staffing), the probability of delay asymptotically equals one, and so the approximate optimization problem reduces to that in an $M/M/1$ system with total capacity $s\mu$ :
\begin{equation}
  \label{eq:mm1}
  \min_s h \frac{\lambda \Gamma(r+1)}{(s\mu-\lambda)^r}+c s.
\end{equation}
The solution to \eqref{eq:mm1} gives us the analog of square-root staffing for this setting:
\begin{equation}
  \label{eq:ss-nl}
  \bar{s}_\lambda=\frac{\lambda}{\mu}+{\sigma}^\star\left( \frac{\lambda}{\mu} \right)^{\frac{1}{r+1}}.
\end{equation}
To establish the result analogous  to the linear case for the performance of this staffing level, we next characterize the  approximation error analogous to \eqref{eq:Q-limit-O(1)}.
\begin{proposition}
\label{prop:accuracy-nl}
For any $\sigma>0$, the following result holds:
\begin{equation}
  \label{eq:nl-moment}
  \lambda \mathbb{E}\xi_\lambda\left( \frac{\lambda}{\mu}+\sigma\left( \frac{\lambda}{\mu} \right)^{\frac{1}{r+1}} \right)=\hat{q}(\sigma)\left(\frac{\lambda}{\mu}\right)^{\frac{1}{r+1}}
+\epsilon_\lambda(\sigma)\left(\frac{\lambda}{\mu}\right)^{\frac{3-r}{2(r+1)}} ,
\end{equation}
where $\hat{q}(\sigma)=\sigma^{-r}\Gamma(r+1)\mu^{1-r}$ and $\epsilon_\lambda(\sigma)=\mathcal{O}(1)$.
\end{proposition}
The proof of this result is postponed to the Appendix. Notice that the approximation error here is $\mathcal{O}(\lambda^{\frac{3-r}{2(1+r)}})$, which decreases in $r$ for $r>1$. For large $r$ values, the many server system operates very similar to the single server system and we find that the approximation gap is small --- for $r=3$, the approximation error is $\mathcal{O}(1)$, which is analogous to the linear case and for $r>3$, the approximation error is in fact $o(1)$.  As $r$ approaches 1, the error in approximating many servers by a single server becomes large and  for $r$ close to $1$, the error is close to $\mathcal{O}(\sqrt{\lambda})$.

Denoting the order of the approximation error in \eqref{eq:nl-moment} by
\begin{equation}
\label{eq:n-lambda}
n_\lambda:=\lambda^{\frac{3-r}{2(1+r)}},
\end{equation}
we now establish the following analog of Theorem~\ref{thm:main} for the convex cost case.
\begin{proposition}
  For convex delay costs with $r>1$,  any staffing level  $\bar{s}_\lambda+\eta_\lambda$, where $\bar{s}_\lambda$ is the staffing prescription defined in \eqref{eq:ss-nl} and $\eta_\lambda=o\left(\sqrt{n_\lambda \lambda^{\frac{1}{1+r}}}\right)$, is $o(n_\lambda)$-optimal for the optimization problem ~\eqref{eq:main-opt-nl}.
\end{proposition}
\begin{proof}
We establish the result by extending the arguments in Proposition~\ref{prop:intuition} and Theorem~\ref{thm:main}. In particular, we note that for $r>1$ and any fixed arrival and service rates, the convexity of the delay cost $\lambda \mathbb{E} \xi_\lambda(s)$ with respect to the number of servers is decreasing. This is easily seen by using \eqref{eq:1} to write
  \begin{align}
    \label{eq:4}
    \lambda \mathbb{E}\xi_\lambda(s)= \frac{\Gamma(r+1)}{(s\mu-\lambda)^{r-1}}\mathbb{E}Q_\lambda(s),
  \end{align}
and then using the established properties of $\mathbb{E}Q_\lambda(s)$. Then, we proceed as in Proposition~\ref{prop:intuition} but define $\Delta_\lambda^Q$ as follows:
\begin{align*}
\Delta_\lambda^Q(i):=\frac{1}{n_\lambda}\left[ \lambda \mathbb{E}\xi_\lambda\left(\frac{\lambda}{\mu}+\sigma \left(\frac{\lambda}{\mu}\right)^{\frac{1}{r+1}}+(i+1) \del\right)-
\lambda\mathbb{E}\xi_\lambda\left(\frac{\lambda}{\mu}+\sigma \left(\frac{\lambda}{\mu}\right)^{\frac{1}{r+1}}+i  \del\right) \right].
\end{align*}
That is,  we replace $\sqrt{\lambda/\mu}$  by $(\lambda/\mu)^{\frac{1}{1+r}}$ and we scale the difference by $n_\lambda$. Then, proceeding as in Proposition~\ref{prop:intuition}, we obtain that for $\del=o\left((n_\lambda \lambda^{\frac{2}{1+r}})^\frac{1}{3}\right)$, we have
\begin{equation}
\label{eq:continuity-nl}
\lim_{\lambda\rightarrow \infty}\left|\epsilon_\lambda(\sigma)-\epsilon_\lambda\left(\sigma+\frac{\del}{(\lambda/\mu)^{\frac{1}{1+r}}}\right)\right|=0.
\end{equation}
Using this we can proceed as in the proof of Theorem~\ref{thm:main} with the appropriate changes to obtain that the staffing level $\bar{s}_\lambda+\eta_\lambda$ is $o(n_\lambda)$-optimal. 
$\qed$
\end{proof}
Thus,  the fact that our approximation  to the delay cost is $\mathcal{O}(n_\lambda)$-accurate, yields that the corresponding prescription has an optimality gap that is an order smaller, i.e., it is $o(n_\lambda)$-optimal. Notice that for $r\ge 3$, we obtain $o(1)$-optimality. 
\subsubsection{Concave delay costs: $r<1$.}
\label{sec:concave} As discussed in \cite{KuR:10}, if $r<1$, the  optimal operating regime is expected to be ``lighter'' than that under the square-root staffing level. However, in such a regime the behavior of the many-server system approaches that of an infinite server queue rapidly, which makes the optimal staffing level only slightly larger than the square-root level. In particular, denoting the optimal staffing level by $s_\lambda^\star$, we have $$\frac{s_\lambda^\star-\lambda/\mu}{(\lambda/\mu)^{\frac{1}{2}+\epsilon}}\rightarrow 0$$ for all $\epsilon>0$. Further, as discussed in \cite{KuR:10}, when dealing with concave delay costs, one expects a last-come-first-serve policy to dominate first-come-first-serve. However, there do not seem to be closed-form expressions for the delay distribution under this policy. So, we do not analyze this case here and leave it for a future study.

\subsection{Application to systems with customer abandonment.}
\label{sec:abandon}
We consider $M/M/s+GI$ systems in which customers have i.i.d. patience times so that a customer whose patience time expires while waiting leaves the system. The  optimization problem with linear holding and capacity costs has been studied in \cite{BaR:10} with a focus on fluid-based prescriptions. We will discuss both fluid- and diffusion-based prescriptions along with their optimality gaps, and we will also discuss solving the case of non-linear delay costs.  For non-linear delay costs, we will restrict attention to the convex case because the observations of Section~\ref{sec:concave}  continue to hold here as well. We  analyze the case of exponential patience times in detail in Section~\ref{sec:exp-pat} and then briefly discuss how these naturally extend to the case of general patience times in Section~\ref{sec:general-pat}.
\subsubsection{Exponential patience times}
\label{sec:exp-pat}
We denote the mean patience time of the customers by $1/\gamma$. Consider the optimization problem \eqref{eq:main-opt} that minimizes the sum of linear holding costs and capacity costs (we ignore abandonment related costs for convenience). 
In this case, as proven in Proposition~5(a) of \cite{BaR:10}, if $h/\gamma>c/\mu$, then the optimal operating regime is critically loaded, i.e., $\frac{s_\lambda^\star}{\lambda/\mu}\rightarrow 1$ as $\lambda \rightarrow \infty$.  The existing results for this model (see, \cite{GMR:02}) suggest that the optimal staffing level for this problem will also have a square-root form. However, there does not appear to be a formal result in the literature that computes the optimal square-root staffing rule that minimizes this cost criterion. (\cite{GMR:02} focuses on staffing to satisfy performance constraints, which is a related problem to the cost minimization that we consider and \cite{ZvZ:12} refines that square-root staffing.) 
Using the results in \cite{GMR:02} and \cite{ZvZ:12}, it is easy to establish that the square-root staffing rule is given by:
$$
\bar{s}_\lambda=\frac{\lambda}{\mu}+\sigma^\star\sqrt{\frac{\lambda}{\mu}},
$$ 
where $\sigma^\star$ solves $$\min_\sigma h q(\sigma)+c \sigma$$ with
$$q(\sigma)=\frac{\mu}{\gamma} \left(\sqrt{\frac{\gamma}{\mu}} H(\sigma\sqrt{\mu/\gamma})-\sigma\right)\left( 1+\sqrt{\frac{\gamma}{\mu}} \frac{H(\sigma\sqrt{\mu/\gamma})}{H(-\sigma)} \right)^{-1},$$
and $H(x)=\frac{\phi(x)}{1-\Phi(x)}$ denotes the hazard rate of the standard normal distribution. Theorem 5 of \cite{ZvZ:12} establishes the corresponding version of Lemma~\ref{lem:mms}, and hence the arguments in Theorem~\ref{thm:main} apply and we obtain that square-root staffing is $o(1)$-optimal in this setting.

Next, we consider non-linear delay costs. The steady-state delay distribution in this system is far more intricate than that in the $M/M/s$ system, which makes computing moments of the steady-state delay difficult. For convenience, we will instead consider delay costs of the form $\left[  \mathbb{E}W(s)\right]^r$ for $r>1$. We expect the insights we obtain, about the order of optimality gaps to be similar, but we leave the exact analysis of delay costs of the form $\mathbb{E}W(s)^r$ for future research. Thus, our non-linear optimization problem is 
\begin{equation}
  \label{eq:nl-a}
  \min_{s\in \mathbb{Z}_+} h \lambda \left[\mathbb{E}W(s)\right]^r+c s.
\end{equation}
Somewhat surprisingly, we find that the optimal solution to \eqref{eq:nl-a} places the system in the overloaded regime with $\frac{s_\lambda^\star}{\lambda/\mu}\rightarrow {\sigma}^\star<1$ for some ${\sigma}^\star\ge 0$. The following result characterizes the optimal fluid prescription and proves that it is $o(1)$-optimal.
\begin{proposition}
For $r>1$,  as $\lambda\rightarrow \infty$, the staffing prescription $\bar{s}_\lambda={\sigma}^\star\frac{\lambda}{\mu}+o(\sqrt{\lambda})$ is $o(1)$-optimal for the optimization problem \eqref{eq:nl-a}, where 
  \begin{equation}
    \label{eq:sigma-pat}
    {\sigma}^\star=\max\left\{ 1-\left( \frac{c \gamma^r}{h r \mu}\right)^{\frac{1}{r-1}} ,0\right\}<1.
  \end{equation}
\end{proposition}
\begin{proof}
Using Little's Law and results in \cite{BaR:10}, for a staffing level $s_\lambda=\sigma \frac{\lambda}{\mu}$, we can write
\begin{equation}
  \label{eq:exp-pat-accuracy}
  \mathbb{E}W_\lambda\left(\sigma\frac{\lambda}{\mu}\right)= \frac{ \mathbb{E}Q_\lambda\left(\sigma\frac{\lambda}{\mu}\right)}{\lambda}=q(\sigma)+\frac{\epsilon_\lambda(\sigma)}{\lambda},~\text{where } \epsilon_\lambda(\sigma)=o(1),
\end{equation}
where $q(\sigma)=\frac{(1-\sigma)^+}{\gamma}$.  Thus, the natural optimization problem that approximates \eqref{eq:nl-a} is
\begin{equation}
  \label{eq:7}
  \min_{0\le \sigma\le 1} h \frac{ (1-\sigma)^r}{\gamma^{r}}+ \frac{c}{\mu} \sigma.
\end{equation}
It is easy to see that \eqref{eq:7} is a convex optimization problem, and further using the first order conditions for optimality, we can verify that $\sigma^\star$ defined in \eqref{eq:sigma-pat} is the unique solution. The results in \cite{BaR:10} immediately imply that the staffing level $\sigma^\star \lambda/\mu$ is $o(1)$-optimal. The arguments in Theorem~\ref{thm:main} can be used to strengthen this result to prove that any staffing level of the form $s_\lambda={\sigma}^\star\frac{\lambda}{\mu}+o(\sqrt{\lambda})$  is $o(1)$-optimal.
  $\qed$
\end{proof}
Notice that for exponential patience times, the approximation error in \eqref{eq:exp-pat-accuracy} itself is $o(1)$, and hence the $o(1)$-optimality directly follows, without the need for further arguments (this is not true for general patience time distributions, where the approximation error is $\mathcal{O}(1)$, but need not be $o(1)$). 

Another point worth noting is that for the case of exponential patience distribution, even though \eqref{eq:exp-pat-accuracy} holds for any $\sigma\ge 0$, for the linear cost structure ($r=1$), the fluid prescription does not give us $o(1)$-optimality. In this case,  if $\frac{h}{\gamma}>\frac{c}{\mu}$, then $\sigma^\star=1$ and the fluid prescription is critically loaded. However, notice that the fluid approximation $q(\sigma)=\frac{(1-\sigma)^+}{\gamma}$ is not continuously differentiable at the $\sigma^\star$, and hence the result breaks down. This observation highlights the importance of differentiability of the approximation  for the asymptotic optimality property.

\subsubsection{General patience times}
\label{sec:general-pat}
For generally distributed patience times, \cite{BaR:10} has shown that under certain technical conditions, the solution to the linear cost minimization problem can lead to an overloaded regime with $\rho>1$. That paper also proves that in that regime, using a fluid approximation, the expected queue-length can be written as
\begin{equation}
\label{eq:ss-fluid}
\mathbb{E}Q_\lambda\left(\sigma \frac{\lambda}{\mu}\right)=q(\sigma)\lambda+\epsilon_\lambda(\sigma),
\end{equation}
where ${q}$ denotes the fluid queue-length and $\epsilon_\lambda=\mathcal{O}(1)$ (unlike the exponential case, the approximation error here does not decrease to zero in general as the arrival rate grows). Thus, the proposed staffing level for such systems is $\bar{s}_\lambda=\sigma^\star \lambda/\mu$, where $\sigma^\star$ solves the fluid optimization problem $\min_\sigma h \mu \bar{q}(\sigma)+c \sigma.$
The analysis in that paper can be easily applied to prove that the relevant continuity property analogous to \eqref{eq:relation} holds so that Theorem~\ref{thm:main} can be used  to obtain the $o(1)$-optimality of the fluid prescription. 

The non-linear cost structure discussed in the exponential case can easily be analyzed for general patience times. The key difference is that the fluid approximation for the expected steady-state delay will be different. Otherwise, the results will be similar and we will obtain that the capacity prescription obtained from the fluid analysis will be $o(1)$-optimal.

\section{Conclusion}
In this paper, we have analyzed optimality gaps of capacity prescriptions for cost minimization in queueing systems, derived using asymptotic analysis. One expects the optimality gap of these prescriptions to be of the same order as that of the error in approximation of the objective function. However,  under some regularity conditions, we find that the optimality gap of such prescriptions is  an order smaller. In particular, we prove that square-root staffing to minimize the sum of linear steady-state holding and capacity costs in $M/M/s$ and $M/M/s+M$ systems is $o(1)$-optimal, i.e., the optimality gap shrinks to zero as the system scale grows without bound. 

The recent paper \cite{Gurvich:14} analyzes a general class of Markovian stochastic systems and establishes rates of convergence for moments of performance measures  to the corresponding diffusion approximations. The results therein imply that diffusion-based expected steady-state queue-length approximations are  $\mathcal{O}(1)$-accurate. If one can establish some regularity conditions along the lines discussed in this paper, for instance that the third derivative is non-negative (or non-positive), then we can further obtain $o(1)$-optimality for these systems. Additional investigation of these regularity conditions would make for a worthwhile future study. We would like to make two additional comments about the work in \cite{Gurvich:14}. First, the main result there applies to more general performance measures, and in particular proves that the approximation error is $\mathcal{O}(1/\sqrt{n})$ for appropriately scaled measures, where $n$ denotes the system scale, and so there is potential to study other performance metrics as well using that result. However,  the result in \cite{Gurvich:14}  does not apply to systems in which the diffusion approximation has a ``reflection'' and consequently, to single server systems (cf., Section~8 of \cite{Gurvich:14}). In this paper, we observed this for the $M/M/s$ system with convex delay costs, in which the optimal operating regime brings the system close to a single-server system, and we identified the approximation error to be $\mathcal{O}(\lambda^{\frac{3-r}{2(1+r)}})$. Thus, the $\mathcal{O}(1)$ approximation result does not hold for $1<r<3$.  The non-linear case in fact  illustrates some differences between systems with and without abandonments. In systems with abandonments, for convex delay costs,  it is optimal to  operate in an overloaded regime. In this case, fluid-based approximations are $\mathcal{O}(1)$-accurate and consequently these prescriptions  are  $o(1)$-optimal. 

Another recent paper that provides a promising direction for future work is \cite{BravermanDai:2015}. That paper considers many-servers systems with customers having phase-type service distribution and exponential patience distribution, i.e., $M/Ph/n+M$ systems. In the Halfin-Whitt asymptotic regime, the authors prove that the approximation error between the appropriately scaled steady-state performance measures computed for the system and for the approximating diffusion is bounded above by $\mathcal{O}(n^{-1/4})$. This result implies that the diffusion-based steady-state queue-length approximations would have an error bounded by $\mathcal{O}(n^{1/4})$. So, if one can establish further regularity conditions as in this paper, we would obtain $o(n^{1/4})$-optimality for these systems. This work is especially encouraging as the framework therein could be potentially applied to other stochastic systems as well.

Finally, we would like to mention that although the paper focuses on unconstrained cost minimization,  our results can be applied to constrained optimization as well. However, we may have situations in which ensuring that the constraints are satisfied may require additional servers than prescribed by square-root staffing and so we may not be able to improve on $\mathcal{O}(1)$-optimality. An example of such a case is in an $M/M/s$ queueing system, for the problem of minimizing the number of servers to ensure that the probability of waiting (before beginning service) is less than a threshold. In this case, one can use the results of \cite{JVZ:11} to show that square-root staffing is only $\mathcal{O}(1)$-optimal.

\bibliographystyle{agsm}      

\bibliography{bibfile}   

\appendix
\section{Proof of Proposition~\ref{prop:mmn-convex}}
For the first part, we will establish that the third derivative of the expected steady-state queue-length with respect to the system utilization is positive for any fixed number of servers. We follow the argument in \cite{Grassman:83}. For convenience, we use $L(\rho)$ to denote the expected number of customers in system in steady-state in the $M/M/s$ system with arrival rate $\lambda$, service rate $\mu$ and utilization $\rho=\frac{\lambda}{\mu s}$. Then, we have
\begin{equation}
\label{eq:erlang}
L=s \rho+ \frac{\rho}{1-\rho} B,
\end{equation}
where $B$ is the probability that all servers are busy.  Following \cite{Grassman:83}, we can write the derivative of $L$ with respect to $\rho$ as
\begin{equation}
\label{eq:dl}
L'=s+(L- s \rho)\left(\frac{s-L+1}{\rho}+\frac{2}{1-\rho}\right).
\end{equation}
Using this relation, we can further differentiate both sides to obtain:
\begin{align}
\label{eq:d2l}
L''&=(L-\ser  \rho ) \left(-\frac{L'}{\rho }+\frac{L-\ser -1}{\rho ^2}+\frac{2}{(1-\rho )^2}\right)-\left(\frac{L-\ser -1}{\rho }-\frac{2}{1-\rho}\right) \left(L'-\ser \right), \text{ and }\\
\label{eq:d3l}
L'''&=-\left(\frac{L-\ser -1}{\rho }-\frac{2}{1-\rho}\right) L''+2 \left(-\frac{L'}{\rho }+\frac{L-\ser -1}{\rho ^2}+\frac{2}{(1-\rho)^2}\right) \left(L'-\ser \right)\\
\nonumber &+(L-\ser  \rho ) \left(-\frac{L''}{\rho }+\frac{2 L'}{\rho ^2}-\frac{2 (L-\ser- 1)}{\rho ^3}+\frac{4}{(1-\rho)^3}\right)
\end{align}
Noting that $\mathbb{E}Q_\lambda=L- \ser  \rho$, we need to prove $L'''>0$ to complete the proof. 

We proceed by substituting the expressions for $L''$ from \eqref{eq:d2l} and that for $L'$ from \eqref{eq:dl}, and finally that for $L$ from \eqref{eq:erlang} into \eqref{eq:d3l}. This yields
\begin{align*}
Z(\ser ,B):=\frac{(1-\rho)^4 \rho^2}{B} L'''&=- \Big(6 \rho ^2 \left(-B^2+(B-1)^3 \rho +3 B-3\right)+(7 B-3) \ser ^2 \rho  (1-\rho)^4\\ 
&-\ser  (1-\rho)^2 (\rho  ((B (12 B-19)+6) \rho -5 B+7)-1)-\ser ^3 (1-\rho)^6\Big).
\end{align*}
To prove $L'''>0$, we need to prove that $Z(\ser ,B)>0$. We will do so by proving $\frac{\partial}{\partial \sen } Z(\sen ,B)>0$ for $1\le \sen  \le \ser$ in Lemma~\ref{lem:partial1} so that $Z(\ser ,B)\ge Z(1,B)$ and then, in Lemma~\ref{lem:partial2}, we will prove that $Z(1,B)>0$ so that we obtain $Z(\ser ,B)>0$. This completes the proof of the first part of the result. We present the proof of the second part of the result after proving these lemmas.

\begin{lemma}
\label{lem:partial1}
We have  $\frac{\partial}{\partial \sen } Z(\sen ,B)>0$ for $1\le \sen  \le \ser$.
\end{lemma}
\begin{proof}
We proceed by computing 
\begin{align*}
f(x):=\frac{1}{(1-\rho)^2}\frac{\partial }{\partial \sen } Z(\sen ,x)= a x^2 - b x+c,
\end{align*}
where
\begin{align*}
a&=12 \rho^2\\
b&=\rho  \Big(14 \sen  (1-\rho)^2+19 \rho +5\Big)\\
c&=3\Big(\sen (1-\rho)^2+\rho\Big)^2+ (3 \rho^2+7 \rho-1).
\end{align*}
We need to prove that $f(B)>0$. Straightforward algebra shows that $c\ge 0$ for $\sen \ge 1$, which gives us $a,b,c\ge 0$ so that $f(x)=0$ has two positive roots. We denote the smaller root by $x_1:=\frac{b-\sqrt{b^2- 4a c}}{2a}$. To prove $f(B)>0$, it suffices to prove that $f'(B)<f'(x_1)$ because $f$ is a convex quadratic function. We will in fact prove that $f'(\bar{B})<f'(x_1)$, where $\bar{B}$ is the following upper bound on $B$:
\begin{align*}
  B \stackrel{(a)}{\le}1+\frac{\ser (1-\rho )^2}{2 \rho }-\frac{(1-\rho ) \sqrt{\ser^2 (1-\rho )^2+4 \ser \rho }}{2 \rho }  \stackrel{(b)}{ \le} 1+\frac{\sen (1-\rho )^2}{2 \rho }-\frac{(1-\rho ) \sqrt{\sen^2 (1-\rho )^2+4 \sen \rho }}{2 \rho }:=\bar{B}.
\end{align*}
The bound $(a)$ follows from the standard Erlang-$C$ bound  (cf. \cite{Harel:88}) and $(b)$ follows by noting that $\sen\le \ser$.
We next compute
\begin{align*}
f'(x_1)&=-\rho(1-\rho)\sqrt{4 \sen  \left(13 \sen  (\rho -1)^2+61 \rho +35\right)+73}\\
f'(\bar{B})&=-\rho(1-\rho)\Big(2 \sen  (1-\rho)+12 \sqrt{\sen  \left(\sen  (\rho -1)^2+4 \rho \right)}+5\Big).
\end{align*}
Noting that $f'(x_1),f'(\bar{B})<0$, we will prove that $\frac{|f'(x_1)|}{\rho(1-\rho)}<\frac{|f'(\bar{B})|}{\rho(1-\rho)}$. Defining $y=\sen (1-\rho)$, we can write
\begin{align*}
g(y,\sen )&:=\left(\frac{f'(\bar{B})}{\rho(1-\rho)}\right)^2-\left(\frac{f'(x_1)}{\rho(1-\rho)}\right)^2\\
&=24 \left(2 y \sqrt{4 \sen +(y-4) y}+5 \sqrt{4 \sen +(y-4) y}+8 \sen +4 y^2-13 y-2\right).
\end{align*}
A straightforward calculation shows that $\frac{\partial g}{\partial \sen }>0$. It follows that for any fixed $y$ and $\sen \ge 1$, $g(y,\sen )$ is minimized  at $\sen =\max\{y,1\}$. So, for the case $y \ge 1$, we can compute
$$
g(y,y)=48(3 y^2-1)>0.
$$
If $y\le 1$, then we use 
$$g(y,\sen )\ge g(y,1)=48(8-y(7-y))>0.
$$
Thus, we have proved that $f(\bar{B})>0$, and because $f$ is a quadratic decreasing function on $[0,\bar{B}]$, it follows that $f(B)>0$ or equivalently $\frac{\partial }{\partial \sen } Z(\sen ,B) >0$ for $1\le \sen  \le \ser$.
$\qed$

\end{proof}

\begin{lemma}
  \label{lem:partial2}
  We have $Z(1,B)>0$.
\end{lemma}
\begin{proof}
  We have
  \begin{align*}
    z(x):=Z(1,x)=\rho & \Big(\rho ^5+\rho  \left(18 x^2+x+12\right)-\rho ^4 (7 x+3)+3 \rho ^3 (x (4 x+3)+3)\\
    &-\rho ^2 (3 x (2 x (x+1)+9)+1)-12 x+6\Big).
  \end{align*}
  Notice that
$$
z''(x)=12 \rho ^2 \left(2 \rho ^2+3-\rho -3 \rho x \right).
$$
It is straightforward to establish that $z''(x)>0$ for $x\le \rho$
(notice that $B\le \rho$ is a well known bound). Further, we have
$$z'(\rho)=-\rho  \Big(12-\rho  (1+\rho  (9-\rho  (\rho +3)))\Big)<0.$$
It follows that $z'(x)<0$ for all $0\le x\le \rho$. Thus, we have
\begin{align*}
  Z(1,B)=z(B)>z(\rho)=6 \rho>0.
\end{align*}
$\qed$\end{proof}

We next prove that for any fixed offered load, the convexity of the expected queue-length with respect to the number of servers $s$ is decreasing. We follow and extend the arguments in \cite{Harel:10}. In particular, we follow the proof of Proposition~5 therein that establishes that the delay probability in an $M/M/s$ queue (that we denote in this paper by $B(\rho,s)$) is convex and decreasing in the number of servers $s$ for any fixed offered load $\lambda/\mu$. It will be convenient to use the terminology in that paper. In particular, we fix $a=\lambda/\mu$ as the offered load and denote the delay probability by $C_s:=B(a/s,s)$. Noting that the expected steady-state queue-length is given by $$
\frac{\lambda}{s\mu-\lambda}B(\rho,s)=\frac{a}{s-a}C_s,
$$
it follows that establishing our result is equivalent to proving that the delay probability $C_s$ has the stated property that its convexity with respect to the number of servers is decreasing. \cite{Harel:10} proves the convexity of the delay probability by establishing that the difference relation $C_{s-1}+C_{s+1}-2C_s>0$ holds. To establish that the convexity is decreasing we prove that
\begin{equation}
\label{eq:third}
(C_{s-1}+C_{s+1}-2C_s) - (C_s+C_{s+2}-2 C_{s+1})=C_{s-1}-3C_s+3C_{s+1}-C_{s+2}>0.
\end{equation}
As in \cite{Harel:10}, we write:
\begin{align*}
C_{s-1}&=\frac{(s-1)(s-a)C_s(a)}{a(s-1-a+C_s)}  \\
C_{s+1}&=\frac{a(s-a)C_s(a)}{s(s+1-a)-a C_s}\\
C_{s+2}&=\frac{a(s+1-a)C_{s+1}}{(s+1)(s+2-a)-aC_{s+1}}.
\end{align*}
Using these relations in \eqref{eq:third} along with the relation $a=\rho s$, we have:
\begin{align}
 C_{s-1}-3C_s+3C_{s+1}-C_{s+2}&=\frac{(1-\rho) \rho ^2 s^2}{2 C_s \rho +s (\rho  (C_s \rho +C_s+1)+(\rho -1) s-3)-2}\nonumber\\
&-\frac{3 (1-\rho) \rho  s}{C_s\rho +(\rho -1) s-1}+\frac{(1-\rho) (s-1)}{\rho  (C_s-\rho  s+s-1)}-3.
  \label{eq:6}
\end{align}
We next determine the sign of each of the denominators of the terms in the above relation. Using the property that $\rho\ge C_s$ and that $0\le C_x\le 1$ for $s-1\le x\le s+2$, we can easily establish that:
\begin{align*}
  2 C_s \rho +s (\rho  (C_s \rho +C_s+1)+(\rho -1) s-3)-2&<0\\
  C_s\rho +(\rho -1) s-1&<0\\
\rho  (C_s-\rho  s+s-1)&>0.
\end{align*}
We next multiply the product of these three terms, which is positive to the terms in \eqref{eq:6}. Thus, to  establish our result, we need to prove that
\begin{align*}
&(C_{s-1}-3C_s+3C_{s+1}-C_{s+2})\times\\
& ( 2 C_s \rho +s (\rho  (C_s \rho +C_s+1)+(\rho -1) s-3)-2)( C_s\rho +(\rho -1) s-1) \rho  (C_s-\rho  s+s-1)>0.
\end{align*}
The term on the left hand side of the above relation can be further simplified to:
\begin{align}
&-(1-\rho)^4  (\rho  (C_s (2 \rho +5)+3)-3)\nonumber\\
& +s (1-C_s\rho) \left(\rho  \left(3 C_s^2 \rho  (\rho +1)-4 C_s (\rho  (4 \rho -5)+4)+3 (2 (\rho -2) \rho +5)\right)-3\right)\nonumber\\
&+(1-\rho)^2 s^2 (\rho  (C_s (\rho  (C_s(5 \rho +7)-8 \rho +9)-16)+2 \rho )+1)\\
&-2 s^3 (1-C_s\rho)^2 ((3 C_s-4) \rho +1)\nonumber\\
&+(1-\rho)^6 s^4. \label{eq:poly}
\end{align}
We next use the bound $\frac{\rho}{s C_s}>\frac{(1-\rho)^2}{(1-C_s)^2}$ (as used in the convexity proof in \cite{Harel:10}) and the additional bound  $C_s>1-(1-\rho)^2\sqrt{\pi s/8}$ (equation (17) in \cite{Harel:10}) that holds for $s\ge 3$, to establish that the term in \eqref{eq:poly} is positive (we omit the details for brevity). This completes the proof.
$\qed$

\section{Proof of Proposition~\ref{prop:accuracy-nl}.}
We first characterize the delay probability and then apply  \eqref{eq:1}.  We use Theorem~1 of \cite{JVZ:11}, which gives us, for any $\rho<1$,
\begin{align}
  \label{eq:prob-bound}
\begin{split}
\left[ \rho+\gamma\left(\frac{\Phi(\alpha)}{\phi(\alpha)}+\frac{2}{3}\frac{1}{\sqrt{s}}+\frac{1}{\phi(\alpha)}\frac{1}{12 s -1}\right)
 \right]^{-1}&\le  B(\rho,s)\\
&\le \left[ \rho+\gamma\left(\frac{\Phi(\alpha)}{\phi(\alpha)}+\frac{2}{3}\frac{1}{\sqrt{s}}\right)
 \right]^{-1},
\end{split}
\end{align}
where 
\begin{align*}
  \alpha&=\sqrt{-2s(1-\rho+\log\rho)}\\
\gamma&=(1-\rho)\sqrt{s}.
\end{align*}

We next focus on the number of servers $s_\lambda=\frac{\lambda}{\mu}+\sigma\left( \frac{\lambda}{\mu} \right)^{\frac{1}{r+1}} $.  We add the subscript $\lambda$ to the terms $\rho,\alpha,\gamma$ to make the dependence on $\lambda$ explicit. Notice that we can rewrite $\gamma_\lambda$ as 
\begin{align}
\label{eq:gamma}
\gamma_\lambda=\sigma \left(\rat\right)^{\frac{1-r}{2(1+r)}}+o\left(\lambda^{\frac{1-r}{2(1+r)}}\right).
\end{align}
Applying Taylor series to $\log \rho_\lambda$, we can express $\alpha_\lambda$ as:
\begin{equation}
  \label{eq:alpha}
\alpha_\lambda=\sigma \left(\rat\right)^{\frac{1-r}{2(1+r)}}+o\left(\lambda^{\frac{1-r}{2(1+r)}}\right).
\end{equation}
Further, noting that $\alpha_\lambda=o(1)$, we can apply Taylor series to $\frac{\Phi(\alpha_\lambda)}{\phi(\alpha_\lambda)}$ to obtain:
\begin{equation}
  \label{eq:2} \frac{\Phi(\alpha_\lambda)}{\phi(\alpha_\lambda)}=\frac{\Phi(0)+\phi(0)\alpha_\lambda+o(\alpha_\lambda)}{\phi(0)+o(\alpha_\lambda)}=\left( \sqrt{\frac{\pi}{2}}+\alpha_\lambda \right)+o(\alpha_\lambda).
\end{equation}
Using \eqref{eq:gamma}, \eqref{eq:alpha} and \eqref{eq:2} in the upper bounding relation in \eqref{eq:prob-bound}, we obtain:
\begin{align}
\label{eq:ub-2}
  B(\rho_\lambda,s_\lambda)\le\left[ \rho_\lambda+\sigma \left(\rat\right)^{\frac{1-r}{2(1+r)}}\left(\left( \sqrt{\frac{\pi}{2}}+\sigma \left(\rat\right)^{\frac{1-r}{2(1+r)}} \right)+\frac{2}{3}\frac{1}{\sqrt{s_\lambda}}\right)\right]^{-1}+o\left(\lambda^{\frac{1-r}{2(1+r)}}\right).
\end{align}
Applying the Taylor series expansion to the first term on the right hand side of \eqref{eq:ub-2} gives us
\begin{align}
\label{eq:ub-3}
  B(\rho_\lambda,s_\lambda)\le 1-\sqrt{\frac{\pi }{2}} \sigma  \left( \rat\right) ^{\frac{1-r}{2(1+r)}}+o\left(\lambda^{\frac{1-r}{2(1+r)}}\right).
\end{align}

Using the same arguments for the lower bounding relation in \eqref{eq:prob-bound}, we obtain:
\begin{align}
\label{eq:lb-2}
  B(\rho_\lambda,s_\lambda)\ge  1-\sqrt{\frac{\pi }{2}} \sigma  \left(\rat\right) ^{\frac{1-r}{2(1+r)}}+o\left(\lambda^{\frac{1-r}{2(1+r)}}\right).
\end{align}
Combining \eqref{eq:ub-3} and \eqref{eq:lb-2}, gives us
\begin{align*}
  B(\rho_\lambda,s_\lambda)= 1-\sqrt{\frac{\pi }{2}} \sigma  \left(\rat\right)^{\frac{1-r}{2(1+r)}}+o\left(\lambda^{\frac{1-r}{2(1+r)}}\right)
\end{align*}

Thus, using \eqref{eq:nl-moment}, we obtain
\begin{equation}
  \label{eq:5}
  \lambda \mathbb{E}\xi_\lambda\left( \frac{\lambda}{\mu}+\sigma\left( \frac{\lambda}{\mu} \right)^{\frac{1}{r+1}} \right)=\hat{q}(\sigma)\left(\frac{\lambda}{\mu}\right)^{\frac{1}{r+1}}
+\epsilon_\lambda(\sigma)\left(\frac{\lambda}{\mu}\right)^{\frac{3-r}{2(r+1)}},
\end{equation}
where  $\epsilon_\lambda(\sigma)=\mathcal{O}(1)$.
$\qed$
\end{document}